\documentclass[12pt]{amsart}
\usepackage{amsmath,amsfonts,amssymb,amsthm}
\usepackage[english]{babel}
\usepackage{latexsym}
\usepackage{amscd, epsfig}
\usepackage{psfrag}
\input epsf.tex

\newcommand{\nn}{\mathbb{N}}

\newcommand{\emp}{\varnothing}
\newcommand{\sq}{\square}

\newcommand{\zz}{\mathbb{Z}}

\newcommand{\rr}{\mathbb{R}}

\newcommand{\pr}{\prime}

\newcommand{\ga}{\gamma}
\newcommand{\si}{\sigma}

\newcommand{\ep}{\epsilon}
\newcommand{\al}{\alpha}
\newcommand{\be}{\beta}

\newcommand{\ve}{\varepsilon}

\newcommand{\vp}{\varphi}

\newcommand{\T}{\mathbf{T}}

\newcommand{\ssu}{\subset}
\newcommand{\sss}{\supset}
\newcommand{\wt}{\widetilde}

\newcommand{\dv}{\divideontimes}

\newcommand{\rt}{{\text {\bf r} } }
\newcommand{\Aut}{{\text {\rm Aut} } }
\newcommand{\A}{{\text {\rm A} } }
\newcommand{\St}{{\text {\rm St} } }
\newcommand{\id}{ {\text{\small{\rm \texttt{I}}}}}
\newcommand{\0}{{\text {\bf 0} } }
\newcommand{\1}{{\text {\bf 1} } }
\newcommand{\io}{\iota}

\newcommand{\Gr}{\mathbb{G}}
\newcommand{\Hr}{\mathbb{H}}
\newcommand{\Br}{\mathbb{B}}

\newcommand{\lan}{\langle}
\newcommand{\ran}{\rangle}

\newcommand{\I}{{\text {\rm i} } }
\newcommand{\V}{{\text {\rm v} } }
\newcommand{\ord}{{\text {\rm ord} } }

\def\ts{\hskip.015cm}

\def\pre{\preccurlyeq}
\def\suc{\succcurlyeq}
\def\any{\ast}

\newtheorem{thm}{Theorem}[section]
\newtheorem{lemma}[thm]{Lemma}

\newtheorem{exc}[thm]{Exercise}
\newtheorem{cor}[thm]{Corollary}
\newtheorem{prop}[thm]{Proposition}
\newtheorem{conj}[thm]{Conjecture}

\newtheorem{Remark}[thm]{Remark}

\begin{document}
\begin{center}
{\large  \bf \sc
Groups of Intermediate Growth:
}
\smallskip

{\large  \bf \sc
an Introduction for Beginners
}

\vskip1.cm

{{\sc Rostislav Grigorchuk}
\smallskip

{\rm
Department of Mathematics \\
Texas A\&M University \\
College Station, TX 77843\\
{\tt \  grigorch@math.tamu.edu}
}

\vskip.8cm

{\sc Igor~Pak}
\smallskip

{\rm Department of Mathematics \\
MIT, Room 2-390\\
Cambridge, MA 02139\\
{\tt \ \ pak@math.mit.edu
}
}}

\vskip.6cm

June 28, 2006
\end{center}
\bigskip

\begin{abstract}
{
We present an accessible introduction to basic results on
groups of intermediate growth.
}
\end{abstract}

\bigskip

\section*{Introduction}

\medskip

\noindent
The study of growth of groups has a long and remarkable
history spanning over much of the twentieth century, and goes
back to Hilbert, Poincare, Alphors, etc.
In~1968 it became apparent that all known classes of groups have
either polynomial or exponential growth, and John Milnor formally
asked whether groups of intermediate growth exist.
The first such examples were
introduced by the first author two decades ago~\cite{G0}
(see also~\cite{G1,G2}), and since then there has been an
explosion
in the number of works on the subject.  While new techniques
and application have been developed, much of the literature
remains rather specialized, accessible only to specialists
in the area.  This paper is an attempt to present the
material in an introductory manner, to the reader
familiar with only basic algebraic concepts.

\smallskip

We concentrate on study of \emph{the first construction},
a finitely generated group~$\Gr$ introduced by the first
author to resolve Milnor's question, and which became
a prototype for further developments.  Our Main Theorem
shows that~$\Gr$ has \emph{intermediate growth}, i.e.
superpolynomial and subexponential.

Our proof is neither the shortest nor gives the best possible
bounds.   Instead, we attempt to simplify the presentation as
much as possible by breaking the proof into a number of
propositions of independent interest, supporting lemmas,
and exercises.    Along the
way we prove two `bonus' theorems: we show that~$\Gr$ is
\emph{periodic} (every element has a finite order) and give
a nearly linear time algorithm for the word problem
in~$\Gr$. We hope that the beginner readers now have an easy
time entering the field and absorbing what is usually viewed
as unfriendly material.

Let us warn the reader that this paper neither gives a survey
nor presents a new proof of the Main Theorem.  We refer to extensive
survey articles~\cite{BGN,BGS,GNS} and a recent book~\cite{Ha} for
further results and references.  Our proofs roughly follow~\cite{G2,MP2},
but the presentation and details are mostly new.

The paper is structured is as follows.
We start with some background information on the growth
of groups (Section~\ref{sec:growth}) and technical results
for bounding the growth function (Section~\ref{sec:tools}).
These technical results have elementary analytic nature; their
proofs are moved to the Appendix (Section~\ref{sec:app}).
In Section~\ref{sec:tree} we study the group~$\Aut(\T)$ of
automorphisms of an infinite binary (rooted) tree.
The `first construction' group~$\Gr$ is introduced in
Section~\ref{sec:grig}, while the remaining
sections~\ref{sec:inf}--\ref{sec:can-lemma} prove the
intermediate growth of~$\Gr$ and one `bonus' theorem.
We conclude with final remarks and few open problems
(Section~\ref{sec:open}).

\smallskip

{\bf Notation.} \ Throughout the paper we use only
{\it left} group multiplication.
For example, a product $\tau_1\cdot\tau_2$ of automorphisms
$\tau_1 , \tau_2 \in \Aut(\T)$ is given by
$[\tau_1\cdot \tau_2](v) = \tau_2(\tau_1(v))$.
We use notation~$g^h = h^{-1} g h$ for conjugate elements,
and~$\id$ for the identity element.
Finally, let $\nn = \{0,1,2,\ldots\}$.

\vskip.8cm

\section{Growth of groups} \label{sec:growth}

Let $S=\{s_1,\ldots,s_k\}$ be a generating set of a
group~$G=\lan S\ran$.  For every group element~$g\in G$,
denote by $\ell(g) = \ell_S(g)$ the length of the
shortest decomposition
$g \, = \, s_{i_1}^{\pm 1} \cdots s_{i_\ell}^{\pm 1}$.  Let
$\ga_G^S(n)$ be the number of elements~$g\in G$
such that~$\ell(g) \le n$. Function~$\ga = \ga_G^S$ is called
the \emph{growth function} of the group~$G$ with respect to
the generating set~$S$. Clearly,
$\ga(n) \, \le \, \sum_{i=0}^n \, (2k)^n \, \le \, (2k+1)^n.$


\begin{exc}  \label{exc:inf}
Let $G$ be an infinite group.  Prove that the growth
function~$\ga$ is
\emph{monotone increasing:} \, $\ga(n+1) \, > \, \ga(n)$, for
all~$n \ge 1$.
\end{exc}


\begin{exc}  \label{exc:subm}
Check that the growth function~$\ga$ is \emph{submultiplicative:}

$\ga(m+n) \, \le \, \ga(m) \, \ga(n)$, for all~$m,n \ge 1$.
\end{exc}


Consider two functions $\ga, \ga^\pr: \nn \to \nn$.
Define $\ga \pre \ga^\pr$ if
$\ga(n) \le C\, \ga^\pr(\al n)$, for all $n>0$ and
some~$C, \al>0$.
We say that~$\ga$ and~$\ga^\pr$ are \emph{equivalent},
write $\ga \sim \ga^\pr$,
if $\ga  \pre \ga^\pr$ and  $\ga^\pr \pre \ga$.


\begin{exc} \label{exc:equiv}
Let $S$ and~$S^\pr$ be two generating sets of~$G$.
Prove that the corresponding growth functions~$\ga_G^S$
and~$\ga_G^{S^\pr}$ are equivalent.
\end{exc}


A function~$f: \nn \to \rr$ is called
\emph{polynomial} if~$f(n) \sim n^\al$, for some~$\al >0$.
A function~$f$ is called \emph{superpolynomial} if
there exists a limit $\frac{\ln \ga(n)}{\ln n} \to \infty$ as
$n \to \infty$.
For example, $n^{\pi}$ is polynomial, while $n^n$ and
$n^{\log\log n}$ are superpolynomial.

Similarly, a function~$f$ is called \emph{exponential} if
$f(n) \sim e^n$. A function~$f$ is called
\emph{subexponential} if there exists a limit
$\frac{\ln \ga(n)}{n} \to 0$ as $n \to \infty$.
For example, $n^ee^n$ and
$\exp\bigl(n/2 - \sqrt{n} \, \log^2 n)$
are exponential, $e^{n/\log n}$ and $n^{\pi}$ are subexponential,
while~$n^n$ is neither.

Let us note also that there are functions which cannot
be categorized.  For example, $\exp\bigl(n^{\sin n}\bigr)$
fluctuates between~1 and $e^n$, so it is neither polynomial nor
superpolynomial, neither exponential nor subexponential.

Finally, a functions~$f$ is said to have
\emph{intermediate growth} if~$f$ is both subexponential
and superpolynomial. For example,
$n^{\log \log n}$, $e^{\sqrt{n}}$, and $e^{n/\log n}$
all have intermediate growth, while
functions $e^{\sqrt{\log n}}$ and
$n!\sim \left(\frac{n}{e}\right)^n \sim e^{n \log n}$
do not.





Exercise~\ref{exc:equiv} implies that
we can speak of groups with \emph{polynomial,}
\emph{exponential} and  \emph{intermediate growth}.
By a slight abuse of notation, we denote by~$\ga_G$ the
growth function with respect to \emph{any} particular
set of generators. Using the equivalence of functions,
we can speak of groups~$G$ and~$H$ as having
\emph{equivalent growth}: $\ga_G \sim \ga_H$.

\begin{exc} \label{exc:power}
Let~$G$ be an infinite group with polynomial growth.
Prove that the direct product~$G^m = G \times G
\times \ldots \times G$ also has polynomial growth,
but $\ga_G \nsim \ga_{G^m}$ for all $m\ge 2$.
Similarly, if~$G$ has exponential growth
then so does~$G^m$, but~$\ga_G \sim \ga_{G^m}$.
\end{exc}

\begin{exc} \label{exc:fi}
Let $H$ be a subgroup of~$G$ of finite index.
Prove that their growth functions are equivalent:
$\ga_H \sim \ga_G$.
\end{exc}

\begin{exc} \label{exc:rate}
Let $S$ be a generating set of a group~$G$, and let
$\ga = \ga_G^S(n)$ be its growth function.
Show that the the limit
$$\lim_{n\to \infty} \frac{\ln \ga(n)}{n}$$
always exist. This limits is called the \emph{growth rate} of~$G$.
Deduce from here that every group~$G$ has either exponential or
subexponential growth.
\end{exc}

\vskip.8cm

\section{Analytic tools} \label{sec:tools}


The following two technical results are key is our analysis of
growth of groups.
Their proofs are based on straightforward analytic arguments and
have no group theoretic content.  For convenience of the reader
we move the proof into Appendix (Section~\ref{sec:app}).


\begin{lemma}[Lower Bound Lemma] \label{lemma:lb}
Let $f: \nn\to \rr_+$ be a monotone increasing function,
such that $f(n) \to \infty$ as $n \to \infty$.
Suppose $f \suc f^m$
for some~$m >1$.  Then $f(n) \suc \exp(n^\al)$ for
some~$\al >0$.
\end{lemma}

For the upper bound, we need to introduce a notation.
Let $f: \nn\to \rr_+$ be a monotone increasing function,
and let:
$$f^{\star \, k}(n) \, := \,
\sum_{(n_1,\ldots,n_k)} \, f(n_1) \cdots f(n_k),
$$
where the summation is over all $k$-tuples
$(n_1,\ldots,n_k) \in \nn^k$ such that $n_1+\ldots+n_k\le n$.

\begin{lemma}[Upper Bound Lemma] \label{lemma:ub}
Let $f(n)$ be a nonnegative monotone increasing function,
such that $f(n) \to \infty$ as $n \to \infty$.
Suppose $f(n) \le C\,f^{\star \, k}(\al \ts n)$ for some~$k \ge 2$,
$C> 0$, and $0 <\al < 1$.  Then $f(n) \pre \exp(n^\be)$ for
some~$\be <1$.
\end{lemma}

Let us note that the functions $f^k$ and $f^{\star \,k}$ are
strongly related:
$$f^k\left(\left\lfloor \frac{n}k\right\rfloor\right) \le
f^{\star \,k}(n) \le n^kf^k(n)$$
However, to analyze the growth we need the lemmas in this
particular form.


\vskip.8cm

\section{Group automorphisms of a tree} \label{sec:tree}

\begin{figure}[hbt]
\begin{center}
\psfrag{0}{$\0$}
\psfrag{1}{$\1$}
\psfrag{00}{$\0\0$}
\psfrag{11}{$\1\1$}
\psfrag{01}{$\0\1$}
\psfrag{10}{$\1\0$}
\psfrag{r}{$\rt$}
\epsfig{file=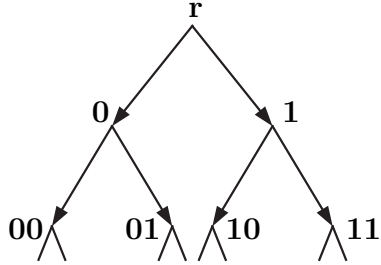,width=5.2cm}
\end{center}
\caption{Infinite binary tree~$\T$.}
\label{f:tree}
\end{figure}

Consider an infinite binary tree~$\T$ as shown in Figure~\ref{f:tree}.
Denote by~$V$ the set of vertices~$v$ in~$\T$, which
are in a natural bijection with finite~$\0$-$\1$ words
$v = (x_0,x_1,\ldots) \in \{\0,\1\}^\ast$.  Note that the root
of~$\T$, denotes~$\rt$, corresponds to an empty word~$\emp$.
Orient all edges in the tree~$\T$ away from the root.
We denote by~$E$ the set of all (oriented) edges in~$\T$.
By definition, $(v,w) \in E$ if $w = v\0$ or $w = v\1$.
Denote by~$|v|$ the distance from the root~$\rt$ to vertex~$v$;
we call it the \emph{level} of~$v$.
Finally, denote by~$\T_v$ a subtree of~$\T$ rooted in~$v \in V$.
Clearly, $\T_v$ is isomorphic to~$\T$.

The main subject of this section is the group $\Aut(\T)$
of automorphisms of~$\T$, i.e. the group of bijections $\tau: V \to V$ which
map edges into edges. Note that the root~$\rt$ is always a fixed
point of~$\tau$. In other words, $\tau(\rt)=\rt$ for all $\tau \in \Aut(\T)$.
More generally, all automorphisms $\tau \in \Aut(\T)$ preserve the level
of vertices:
$|\tau(v)| = |v|$, for all~$v\in V$. Denote by~$\id \in\Aut(\T)$
a trivial (identity) automorphism of~$\T$.

\begin{figure}[hbt]
\begin{center}
\psfrag{0}{$\0$}
\psfrag{1}{$\1$}
\psfrag{T0}{$\T_\0$}
\psfrag{T1}{$\T_\1$}
\psfrag{r}{$\rt$}
\psfrag{a}{$a$}
\epsfig{file=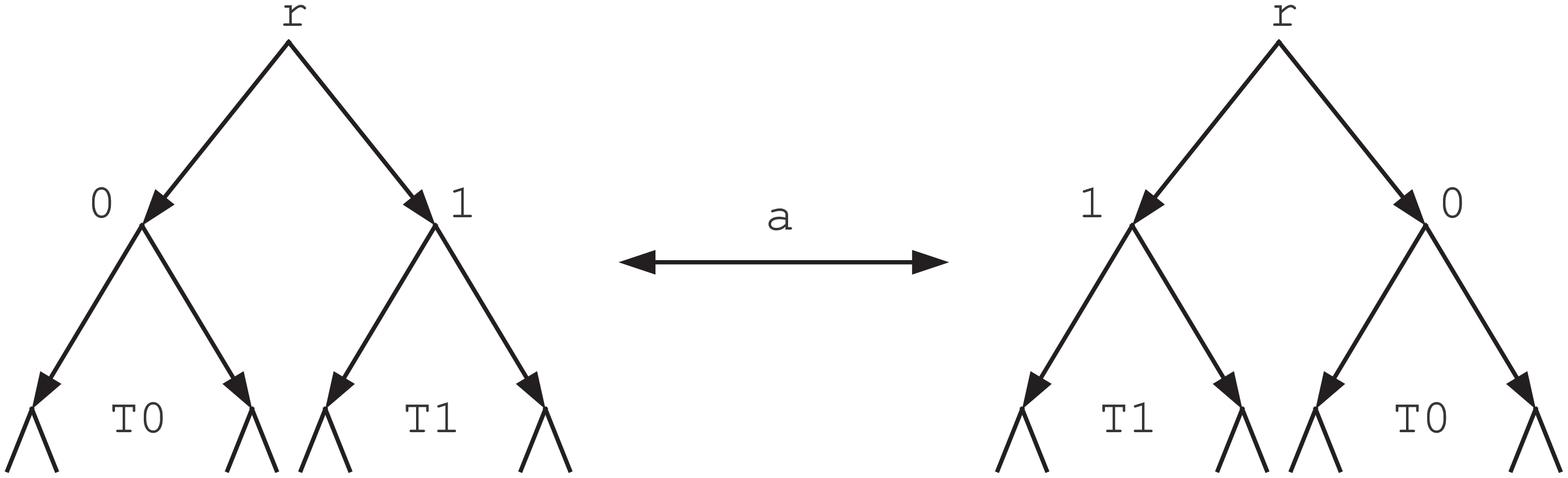,width=11.2cm}
\end{center}
\caption{Automorphism~$a \in \Aut(\T)$.}
\label{f:aut-a}
\end{figure}

An example of a nontrivial automorphism $a \in \Aut(\T)$ is given
in Figure~\ref{f:aut-a}. This is the most basic automorphism which
will be used throughout the paper, and can be formally defined
as follows. Denote by $\T_\0$ and $\T_\1$ the left and right
subtrees (branches) of the tree~$\T$, with roots at~$\0$ and~$\1$,
respectively. Let~$a$ be an automorphism which maps~$\T_\0$
into~$\T_\1$ and preserves the natural order on vertices:
$$
\tau: \, (\0,x_1,x_2,\ldots) \longleftrightarrow
(\1,x_1,x_2,\ldots).
$$
Clearly, automorphism~$a$ is an involution: $a^2=\id$.

Similarly, one can define an automorphism $a_v$ which
exchanges two branches $\T_{v\0}$ and $\T_{v\1}$ of a subtree
$\T_v$ rooted in~$v \in V$.  These automorphisms will be used
in the next section to define finitely generated subgroups
of~$\Aut(\T)$.

More generally, denote by $\Aut(\T_v)$ the subgroup of automorphisms
in~$\Aut(\T)$ which preserve subtree~$T_v$ and are trivial on
the outside of~$T_v$.  There is a natural graph isomorphism
$\io_v : \T \to \T_v$ which extends to a group isomorphism
$\io_v: \Aut(\T)\to \Aut(\T_v)$.

By definition, every automorphism $\tau \in \Aut(\T)$
maps two edges leaving vertex~$v$ into two edges leaving vertex~$\tau(v)$.
Thus we can define a \emph{sign}
$\ep_v(\tau) \in \{0,1\}$ as follows:
$$
\ep_v(\tau) = \left\{ \aligned
& 0, \ \ \ \text{if} \ \ \, \tau(v\0) = \tau(v)\0, \ \, \tau(v\1) = \tau(v)\1,
\\
& 1, \ \ \ \text{if} \ \ \, \tau(v\0) = \tau(v)\1, \ \, \tau(v\1) = \tau(v)\0.
\endaligned \right.
$$
In other words, $\ep_v(\tau)$ is equal to~$0$ if the automorphism maps
the left edge leaving vertex~$v$ into the left edge leaving~$\tau(v)$,
and is equal to~$1$ if the automorphism maps
the left edge leaving~$v$ into the right edge leaving~$\tau(v)$.

Observe that the collection of signs $\{\ep_v(\tau), v\in \T\}$ can
take all possible values, and uniquely determines the
automorphism~$\tau\in \Aut(\T)$. As a corollary, the
group~$\Aut(\T)$ is uncountable and cannot be finitely
generated.

\smallskip

\begin{figure}[hbt]
\begin{center}
\psfrag{0}{$\0$}
\psfrag{1}{$\1$}
\psfrag{r}{$\rt$}
\psfrag{0}{$\0$}
\psfrag{t2}{$\tau_1$}
\psfrag{t1}{$\tau_0$}
\epsfig{file=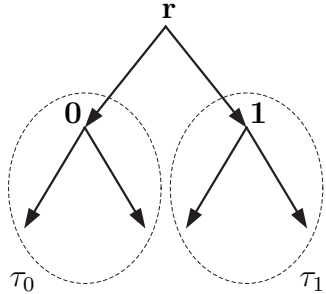,width=6.2cm}
\end{center}
\caption{Automorphism~$\tau = \vp(\tau_0,\tau_1) \in \Aut(\T)$.}
\label{f:psi}
\end{figure}

To further understand the structure of~$\Aut(\T)$, consider a map
$$
\vp: \ \Aut(\T) \times \Aut(\T) \, \to \, \Aut(\T),
$$
defined as follows.  If $\tau_0, \tau_1 \in \Aut(\T)$, let
$\tau = \vp(\tau_0,\tau_1)$ be an automorphism defined by
$\tau := \io_\0(\tau_0) \cdot \io_\1(\tau_1) \in \Aut(\T)$.
Here $\io_\0(\tau_0) \in \Aut(\T_\0)$ and
$\io_\1(\tau_1) \in \Aut(\T_\1)$ are the automorphisms of
subtrees~$\T_\0$ and~$\T_1$, respectively, defined as above.
Pictorially, automorphism~$\tau$ is shown in Figure~\ref{f:psi}.

For any group~$G$, the \emph{wreath product} $G \wr \zz_2$ is
defined as a semidirect product $(G \times G) \rtimes \zz_2$,
with~$\zz_2$ acting by exchanging two copies of~$G$.


\begin{prop} $\Aut(\T) \simeq \Aut(\T) \wr \zz_2$.
\end{prop}


\begin{proof} Let us extend the map $\vp$ to an isomorphism
$$
\vp: \ \bigl(\Aut(\T) \times \Aut(\T)\bigr) \rtimes \zz_2
\, \longrightarrow \, \Aut(\T)
$$
as follows.  When $\si =\id$, let
$\vp(\tau_0,\tau_1;\si) := \vp(\tau_0,\tau_1)$, as before.
When $\si \ne \id$, let
$\vp(\tau_0,\tau_1;\si) := \vp(\tau_0,\tau_1)\cdot a$,
where~$a \in \Aut(\T)$ is defined as
above.  Now check that multiplication of automorphisms
$\vp(\cdot)$ coincides with that of the semidirect product,
and defines the group isomorphism.
We leave this easy verification to the reader.
\end{proof}


We denote by $\psi = \vp^{-1}$ the
isomorphism $\psi:  \Aut(\T) \to \Aut(\T)\wr \zz_2$ defined
in the proof above.
This notation will be used throughout the paper.

\smallskip

\begin{exc} \label{exc:autm} Let $\A_m \ssu \Aut(\T)$ be a
subgroup of all automorphisms~$\tau \in \Aut(\T)$ such that
$\ep_v(\tau) = 0$ for all~$|v|\ge m$. For example, $\A_1 = \{\id,a\}$.
Use the idea above to
show that
$$\A_m \simeq \zz_2 \, \wr \, \zz_2 \, \wr \, \cdots \, \wr \, \zz_2 \ \ \
(m \ \text{times}).$$
Conclude from here that the order of~$\, \A_m$ is
$\, |\,\A_m| \, = \, 2^{2^m-1}$.
\end{exc}

\begin{exc} \label{exc:inf-el} Consider a (unique)
tree automorphism~$\tau \in \Aut(\T)$ with signs
given by:
$\ep_v(\tau)=1$ if $v = \1^k = \1\ldots\1$~($k$ times),
for~$k \ge 0$, and $\ep_v(\tau) = 0$ otherwise.
Check that~$\tau$ has infinite order in~$\Aut(\T)$.

\noindent
{\bf Hint:}  {\rm Consider elements $\tau_m \in \A_m$ with
signs as in the definition of above, and~$k <m$.
Show that the order~$\ord(\tau_m) \to \infty$ as~$m \to \infty$,
and deduce the result from here.}
\end{exc}

\vskip.8cm

\section{The first construction} \label{sec:grig}

In this section we define a finitely generated
group~$\Gr \ssu \Aut(\T)$ which we call \emph{the
first construction}.  Historically, this is the
first example of a group with  intermediate growth~\cite{G0}.

Let us first define group~$\Gr$ implicitly, by
specifying the necessary conditions on generators.
Let $\Gr = \lan a,b,c,d \ran \ssu \Aut(\T)$, where~$a$
is the automorphism defined as in Section~\ref{sec:tree},
and automorphisms~$b$,~$c$,~$d$ satisfy the
following conditions:
$$
(\circ) \ \ \ b = \vp(a,c), \ \ c = \vp(a,d),
\ \ d = \vp(\id,b).
$$
Observe that the automorphisms~$b$,~$c$, and~$d$ are defined
through each other. Since the generator~$d$ is acting as identity
automorphism on the left subtree~$\T_\0$, and as~$b$ on the right
subtree~$\T_\1$, one can recursively compute the action of all
three automorphisms~$b,c,d \in \Aut(\T)$.

Here is a direct way to define automorphisms $b,c,d$~:
$$(\ast) \ \ \ \
\aligned
b \, & :=  \,  (a_{\0} \ \cdot  \, a_{\1^3\0} \,
\cdot \, a_{\1^6\0} \cdot \ldots) \,
(a_{\1\0} \ \cdot  \, a_{\1^4\0} \,
\cdot \, a_{\1^7\0} \cdot \ldots ),
\\
c \, & := \, (a_{\0} \, \cdot \, a_{\1^3\0} \,
\cdot \, a_{\1^6\0} \cdot \ldots ) \,
(a_{\1^2\0} \ \cdot  \, a_{\1^5\0} \,
\cdot \, a_{\1^8\0} \cdot \ldots) ,
\\
d \, & := \,(a_{\1\0} \ \cdot  \, a_{\1^4\0} \,
\cdot \, a_{\1^7\0} \cdot \ldots )
(a_{\1^2\0} \, \cdot \, a_{\1^5\0} \,
\cdot \, a_{\1^8\0} \cdot \ldots ),
\endaligned
$$
where~$\1^m$ is short for $\1\ldots\1$ ($m$~times).
Note that the automorphisms $a_{\1^m\0}$ used
in~$(\ast)$ commute with each other, and thus
elements~$b,c,d \in \Aut(\T)$ are well defined.

\begin{figure}[hbt]
\begin{center}
\psfrag{b}{$b$}
\psfrag{c}{$c$}
\psfrag{d}{$d$}
\epsfig{file=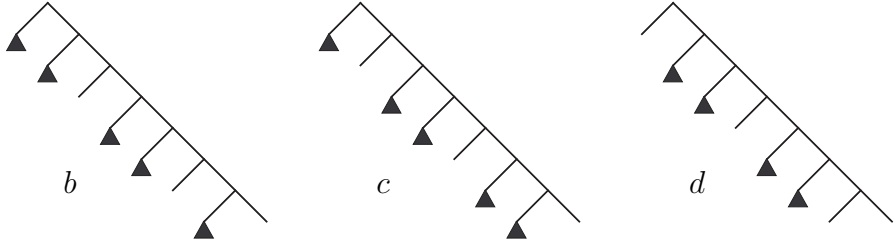,width=12.cm}
\end{center}
\caption{Elements $b,c$ and~$d \in \Aut(\T)$.}
\label{f:bcd}
\end{figure}

Elements~$b,c,d \in \Aut(\T)$ are graphically shown
in Figure~\ref{f:bcd}.  Here black triangles in vertices
of trees represent the subtrees swaps.


\begin{thm}[Main Theorem] \label{thm:main}
Group $\Gr = \lan a,b,c,d\ran$ has intermediate
growth.
\end{thm}


The proof of Theorem~\ref{thm:main} is quite involved and
occupies much of the rest of the paper.


\begin{exc} \label{exc:circ}
Check that elements $b,c,d \in \Aut(\T)$
defined by~$(\any)$ satisfy conditions~$(\circ)$.
\end{exc}


\begin{exc} \label{exc:bcd}
Check that elements $b,c$, and~$d$ are involutions (have order~2),
commute with each other, and satisfy \, $b\cdot c\cdot d=\id$.
Conclude from here that $\lan b,c,d\ran \simeq \zz_2^2$ and that the
group $\Gr = \lan a,b,c,d \ran$  is $3$-generated.
\end{exc}


\begin{exc} \label{exc:order}
Check the following relations in~$\Gr :$ \,
$(ad)^4 = (ac)^8 = (ab)^{16}  = \id$.
Deduce from here that 2-generator subgroups
$\lan a, b\ran,\lan a,c\ran,\lan a,d\ran \ssu \Gr$ are finite.
\end{exc}

While these exercises have straightforward `verification style'
proofs, they will prove useful in the future. Thus
we suggest the reader studies them before proceeding to read
(hopefully) the rest of the paper.


\vskip.8cm

\section{Group $\Gr$ is infinite} \label{sec:inf}

We have yet to establish that~$\Gr$ is infinite.
Although one can prove this directly, the proof below
introduces definitions and notation which will be helpful
in the future.

\smallskip

Let $\St_\Gr(n)$ denote a subgroup of~$\Gr$ stabilizing
all vertices with level~$n$.  In other words, $\St_\Gr(n)$
consists of all automorphisms~$\tau \in \Gr$ such that
$\tau(v)=v$ for all vertices~$v \in \T$ with~$|v| = n$:
$$
\St_\Gr(n) \, = \, \bigcap_{|v| = n} \St_\Gr(v).
$$
The subgroup $\Hr := \St_\Gr(1)$ is called the
\emph{fundamental subgroup} of~$\Gr$.


\begin{lemma} \label{lemma:hgen}
Let~$\Hr \ssu \Gr$ be the fundamental subgroup defined above.
Then:
$$\Hr = \lan b,c,d,b^a,c^a,d^a \ran, \quad \Hr \lhd \Gr,
\ \ \ \text{and} \ \ \  [\Gr:\Hr]=2.$$
\end{lemma}

\begin{proof} From Exercise~\ref{exc:bcd} we conclude that
every reduced decomposition~$w$ is a product
$w = (a) \any a \any a \any \ldots \any a \any (a)$, where
each~$\any$ is either~$b$,~$c$, or~$d$, while the first
and last~$a$ may or may not appear.  Denote by $|w|$ the
length of the word~$w$, and by~$|w|_a$ the number of
occurrences of~$a$ in~$w$.  Note that $w \in \Hr$ if and only
if~$|w|_a$ is even. This immediately implies the third part
of the lemma.  Since every subgroup of index~2 is normal
this also implies the second part.

For the first part, suppose $|w|_a$ is even. Join subsequent
occurrences of~$a$ to obtain~$w$ as a product of~$\any$
and~$(a \any a)$.  Since~$a^2=\id$, we have~$(a\any a) = \any^a$,
which implies the result.
\end{proof}

This following exercise generalizes the second part of
Lemma~\ref{lemma:hgen} and will be used in
Section~\ref{sec:upper}
to prove the upper bound on the growth function of~$\Gr$.

\begin{exc} \label{exc:stab-ind}
Check that the stabilizer subgroup $\Hr_n:=\St_\Gr(n)$
has finite index in~$\Gr:$
$[\Gr:\Hr_n] \le |\,\A_n | = 2^{2^n-1}$ {\rm(see Exercise~\ref{exc:autm}).}
\end{exc}

\smallskip


Let~$\psi = \vp^{-1}: \Aut(\T) \to
\bigl(\Aut(\T) \times \Aut(\T)\bigr)\rtimes \zz_2$
be the isomorphism defined in Section~\ref{sec:tree}.
By definition, $\Hr \ssu \Gr \ssu \Aut(\T)$.


\begin{lemma} \label{lemma:surj}
The image $\psi(\Hr)$ is a subgroup
of~$\Gr \times \Gr$, such that projection of $\psi(\Hr)$
onto each component is surjective.
\end{lemma}

\begin{proof} By definition, $\Hr$ stabilizes~$\0$ and~$\1$,
so $\psi(\Hr) \ssu \Aut(\T) \times \Aut(\T)$.
From Exercise~\ref{exc:circ} we have
$$
\psi : \ \ \left\{ \aligned
& \ b \to (a,c) \,, \ \qquad b^a \to (c,a) \,, \\
& \ c \to (a,d) \,, \ \qquad c^a \to (d,a) \,,\\
& \ d \to (\id,b) \,, \ \qquad d^a \to (b,\id) \,.
\endaligned \right.
$$
Now Lemma~\ref{lemma:hgen} implies that
$\psi(\Hr) \ssu \Gr \times \Gr$.
On the other hand, the projection of~$\psi(\Hr)$ onto each component
contains all four generators~$a,b,c,d \in \Gr$, and is therefore
surjective.
\end{proof}


\begin{cor} \label{cor:inf}
Group~$\Gr$ is infinite.
\end{cor}

\begin{proof} From Lemma~\ref{lemma:hgen} and
Lemma~\ref{lemma:surj} above, we have $\Hr$ is a
proper subgroup of~$\Gr$ which is mapped surjectively
onto~$\Gr$.  If~$|\Gr| <\infty$, then
$|\Gr| > |\Hr| \ge |\Gr|$, a contradiction.
\end{proof}

\medskip

Here is a different application of Lemma~\ref{lemma:surj}.
Let $G \ssu \Aut(\T)$ be a subgroup of the group
automorphisms of the binary tree~$\T$.
Denote by $G_v = \St_G(v)|_{\T_v} \ssu \Aut(\T_v)$
the subgroup of~$G$ of elements which fix vertex~$v \in \T$
with the action restricted only to the subtree~$\T_v$.  We say
that~$G$ has (strong) \emph{self-similarity property}
if $G_v \simeq G$ for all~$v \in \T$.

\begin{cor} \label{cor:self}
Group~$\Gr$ has self-similarity property.
\end{cor}

\begin{proof}
Use the induction on the level~$|v|$.
By definition, $\Gr_\rt = \Gr$, and by Lemma~\ref{lemma:surj}
we have $\Gr_\0, \Gr_1 \simeq \Gr$. For general~$v \in \T$
we similarly have $\Gr_{v\0}, \Gr_{v\1} \simeq \Gr_v$.
This implies the result.
\end{proof}

\begin{exc} \label{exc:infinite}
Consider the following rewriting rules:
$$\eta: \, a \to aba, \ \ b \to d, \ \ c \to b, \ \ d \to c$$
Define a sequence of elements in~$\Gr$: $x_1 = a$ and
$x_{i+1} := \eta(x_i)$
for all $i \ge 1$.  Prove directly that all these elements are
distinct.  Conclude from here that~$\Gr$ is infinite.
\end{exc}

\vskip.8cm

\section{Superpolynomial growth of~$\Gr$} \label{sec:lower}

In this section we prove the first half of Theorem~\ref{thm:main},
by showing that the growth function~$\ga$ of group~$\Gr$ satisfies
conditions of the Lower Bound Lemma.

\smallskip

Two groups $G_1$ and $G_2$ are called \emph{commeasurable},
denotes $G_1 \approx G_2$, if they contain isomorphic
subgroups of finite index:
$$H_1 \ssu G_1, \, H_2 \ssu G_2, \, \ H_1 \simeq H_2, \ \
\text{and} \ \ [G_1 : H_1], \, [G_2 : H_2] < \infty.$$
For example, group~$\zz$ is commeasurable with the infinite
dihedral group~$\textsc{D}_\infty \simeq \zz \rtimes \zz_2$.
Of course, all finite groups are commeasurable to each
other. Another example is $\Hr \approx \Gr$, since~$\Hr$ is
a subgroup of finite index in~$\Gr$.
Note also that commeasurability is an equivalence
relation.


\begin{prop} \label{prop:ss}
Groups~$\Gr$ and~$\Gr \times \Gr$ are commeasurable:
$\Gr \approx \Gr \times \Gr$.
\end{prop}

Proposition~\ref{prop:ss} describes an important phenomenon
which can be formalized as follows.   The group~$G$ is called
\emph{multilateral} if~$G$ is infinite and~$G \approx G^m$
for some~$m \ge 2$.  As we show below, all such groups have
superpolynomial growth.

\smallskip

To prove the proposition, consider the subgroups~$\Hr \ssu \Gr$
and~$\wt{\Hr} := \psi(\Hr) \ssu \Gr \times \Gr$.
By Lemma~\ref{lemma:hgen} we have~$[\Gr:\Hr] < \infty$.
Since~$\psi$ is a group isomorphism,
we also have~$\wt \Hr \simeq \Hr$.
If we show that $[\Gr\times \Gr:\wt \Hr] < \infty$,
then $\Gr \approx \Gr \times \Gr$, as claimed in
Proposition~\ref{prop:ss}.

\smallskip

Denote by~$\Br=\lan b \ran^G$ the normal closure
of~$b \in\Gr$, defined
as $\Br := \lan g^{-1} b g \mid g \in \Gr\ran$.


\begin{lemma} \label{lemma:Br}
Subgroup~$\Br$ has a finite index in~$\Gr$.  More precisely,
$[\Gr:\Br] \le 8$.
\end{lemma}

\begin{proof}
By Exercise~\ref{exc:order}, we have $a^2 = d^2 = (ad)^4 = \id$.
It is easy to see now that the 2-generated subgroup
$\lan a,d\ran \ssu \Gr$ is a dihedral group~$D_4$ of order~8.
By Exercise~\ref{exc:bcd}, we have $\Gr = \lan a,b,d\ran$.
Therefore, $\Gr/\Br$ is a quotient of~$\lan a,d\ran$, and
$[\Gr : \Br] \le |D_4| = 8$.
\end{proof}


\begin{lemma} \label{lemma:BrBr}
Subgroup $\wt\Hr = \psi(\Hr)$
contains~$\Br \times \Br \ssu \Gr \times \Gr$.
\end{lemma}

\begin{proof}
By Lemma~\ref{lemma:hgen}, we know that
$\wt \Hr \sss \lan\psi(d), \psi(d^a) \ran = \lan (1,b), (b,1)\ran$.
Let $x \in \Hr$ and $\psi(x) = (x_0,x_1)$.  We have:
$$\aligned
\psi(d^x) \, & = \, \psi(x^{-1} d x) \, = \,
\psi(x^{-1}) \, \psi(d) \, \psi(x) \, = \,
(x_0^{-1},x_1^{-1}) \, (\id,b) \, (x_0,x_1) \\
 \, & = \,
(\id,x_1^{-1}b x_1) \, = \, (\id,b^{x_1}).
\endaligned
$$
By Lemma~\ref{lemma:surj}, here we can take any element~$x_1 \in \Gr$.
Therefore, the image~$\psi(\Hr)$ contains all elements of the
form~$(\id,b^g)$, $g \in \Gr$. By definition, these elements generate
a subgroup~$1 \times \Br$. In other words,
$\wt \Hr = \psi(\Hr) \sss 1 \times \Br$.  Similarly,
using the element~$d^a$ in place of~$d$, we obtain
$\wt \Hr \sss \Br \times 1$.  Therefore, $\wt \Hr \sss \Br \times \Br$,
as desired.
\end{proof}


Now Proposition~\ref{prop:ss} follows immediately once we note
that $\Br \times \Br \ssu \wt \Hr \ssu \Gr \times \Gr$, and
by Lemma~\ref{lemma:Br} the index
$$
[\Gr\times\Gr : \wt \Hr] \, \le \, [\Gr\times\Gr : \Br \times \Br]
\, = \, [\Gr : \Br]^2 \, \le \, 8^2 \, =  \, 64.
$$
Since~$\Gr$ is infinite (Corollary~\ref{cor:inf}) this
implies that group~$\Gr$ is multilateral. \ $\sq$

\smallskip

\begin{lemma} \label{lemma:lower}
Every multilateral group~$G$ has superpolynomial growth.
Moreover, the growth function $\ga_G(n) \suc \exp(n^\al)$
for some~$\al > 0$.
\end{lemma}

\begin{proof}
By definition, $G$ is infinite, and $G \approx G^m$ for some $m >1$.
In other words, there exist $H \ssu G$, $\wt H \ssu G^m$ such that
$H \simeq \wt H$ and $[G:H], [G^m:\wt H]<\infty$.
From Exercise~\ref{exc:fi} we obtain
$\ga_G \sim \ga_{H} \sim \ga_{\wt H} \sim \ga_{G^m}$.
Thus~$\ga_G \suc \ga_{G^m}$, and the Lower Bound
(Lemma~\ref{lemma:lb}) implies the result.
\end{proof}

\smallskip

Now Proposition~\ref{prop:ss} and Lemma~\ref{lemma:lower} immediately
imply the first part of Theorem~\ref{thm:main}:

\begin{cor} \label{cor:lower}
Group~$\Gr$ has superpolynomial growth.
Moreover, the growth function $\ga_{\Gr}(n) \suc \exp(n^\al)$
for some~$\al > 0$.
\end{cor}

\vskip.8cm

\section{Length of elements and rewriting rules} \label{sec:length}

To prove the second half of Theorem~\ref{thm:main}
we derive sharp upper bounds on the
growth function~$\ga=\ga_\Gr^S$ of the group~$\Gr$ with
the generating set~$S = \{a,b,c,d\}$.  In this section
we obtain some recursive bounds on the
length~$\ell(g) = \ell_\Gr^S(g)$ of elements~$g \in \Gr$
in terms of~$S$.  Note that
although~$\Gr$ is 3-generated, having the fourth generator
is convenient for technical reasons.

\smallskip

We begin with a simple classification of reduced
decompositions of elements of~$\Gr$ following the approach
in the proof of Lemma~\ref{lemma:hgen}.  We define
four \emph{types} of reduced decompositions:

\smallskip

$(\I)$ \ \  $\ $ if \ $g = a \any a \any a \cdots \any a \any a$,

$(\I\I)$  \  $\ $ if \ $g = a \any a \any a \cdots \any a \, \any$,

$(\I\I\I)$  \ \ if \ $g = \any \, a \any a \any \cdots a \any a$,

$(\I\V)$  \ \ if \ $g = \any \, a \any a \any \cdots a \any a \, \ts \any$.

\smallskip

\noindent
Of course, element~$g$ can have many different reduced decompositions.
On the other hand, the type of a decomposition is almost completely
determined by~$g$.

\begin{lemma}\label{lemma:type}
Every group element $g \in \Gr$ has all of its reduced
decompositions of the same type~$(\I)$, \emph{or} of type~$(\I\V)$,
\emph{or} of type~$(\I\I)$ and~$(\I\I\I)$.
\end{lemma}

\begin{proof} Recall that the number of $a$'s in a reduced
decomposition of~$g \in \Gr$ is even if $g \in \Hr$, and is
odd otherwise. Thus~$g$ cannot have decompositions of
type~$(\I)$ and $(\I\V)$ at the same time.  Noting that
decompositions of type~$(\I)$ and~$(\I\V)$ have odd length
while those of type~$(\I\I)$ and~$(\I\I\I)$ have even length
implies the result.
\end{proof}

It is easy to see that one cannot strengthen Lemma~\ref{lemma:type}
since some elements can have decompositions of both
type~$(\I\I)$ and~$(\I\I\I)$.  For example, $adad = dada$ by
Exercise~\ref{exc:bcd}, and both are reduced decompositions.
From this point on we refer to elements $g\in \Gr$ as of \emph{type}
$(\I)$, $(\I\I/\I\I\I)$, or~$(\I\V)$ depending on the type of their
reduced decompositions.

\smallskip

In the next lemma we use the isomorphism
$\psi = \vp^{-1} : \Aut(\T) \to \Aut(\T) \wr S_2$,
where $S_2 = \{\id,a\} \simeq \zz_2$.

\begin{lemma} \label{lemma:half}
Let $\ell(g)$ be the length of $g \in \Gr$ in generators~$S=\{a,b,c,d\}$.
Suppose $\psi(g) = (g_0,g_1; \si)$, where $g_0,g_1 \in \Gr$
and~$\si \in S_2$. Then:

\smallskip

$\ell(g_0),\ell(g_1) \le \frac12(\ell(g)-1)$ if $g$ has type~$(\I)$,

$\ell(g_0),\ell(g_1) \le \frac12\,\ell(g)$ if $g$ has type~$(\I\I/\I\I\I)$,

$\ell(g_0),\ell(g_1) \le \frac12(\ell(g)+1)$ if $g$ has type~$(\I\V)$.
\end{lemma}

\begin{proof}
Fix an element $g \in \Gr$, and let $g_0,g_1,\si$ be as in the lemma.
We have $\si = \id$ if $g \in \Hr$, and $\si = a$ otherwise (see
the proof of Lemma~\ref{lemma:hgen}).  For every reduced decomposition
$w = (a) \any a \any a \cdots \any a \any (a)$ of~$g$
we shall construct decompositions of elements~$g_0,g_1$ with lengths
as in the lemma.
As before, we use~$\any$ to denote either of the
generators~$b,c,d$.  Also, for every~$\any$ in a reduced decomposition
denote by $\pi(\any)$ the number of~$a$'s preceding~$\any$.

Consider the following \emph{rewriting rules}:
$$\Phi_0: \ \qquad \left\{
\aligned
& a \to \id, \\
& b \to a, \ \  c \to a, \ \ d \to \id \ \ \ \text{if} \ \ \, \pi(\any)
\ \ \,
\text{is odd,}\\
& b \to c, \ \  c \to d, \ \ \ts \ts d \to b \ \ \ \text{if} \ \ \, \pi(\any)
\ \ \,
\text{is even,}
\endaligned\right.
$$
and
$$\Phi_1: \  \qquad \left\{
\aligned
& a \to \id, \\
& b \to a, \ \  c \to a, \ \ d \to \id \ \ \ \text{if} \ \ \, \pi(\any)
\ \ \,
\text{is even,}\\
& b \to c, \ \  c \to d, \ \ \ts \ts \ts d \to b \ \ \ \text{if} \ \ \, \pi(\any)
\ \ \,
\text{is odd.}
\endaligned\right.
$$
These rules act on words~$w$ in generators~$S$, and substitute each
occurrence of a letter with the corresponding letter or~$\id$.

Let $\Phi_0(w)$, $\Phi_1(w)$ be the words obtained from the word
$w = (a) \any a \cdots a \any (a)$ by the rewriting rules
as above, and let $g^\pr_0, g^\pr_1 \in \Gr$ be group elements defined
by these products.  Check by induction on the length~$\ell(g)$
that~$\psi(g) = (g^\pr_0,g^\pr_1;\si)$.  Indeed, note that the rules
give the first and second components in the formula for~$\psi$
in the proof of Lemma~\ref{lemma:surj}.  Now, as in the proof of
Lemma~\ref{lemma:hgen} subdivide the product~$w$ into elements~$(a)$
and~$(\any \, \ts a \ts \any)$, and obtain the induction step.
From here we have $g_0 = g^\pr_0$, $g_1 = g^\pr_1$,
and by construction of rewriting rules the lengths of $g_0,g_1$
are as in the lemma.
\end{proof}

As we show below, the rewriting rules are very useful
in the study of group~$\Gr$, but also in a more general
setting.

\begin{cor}\label{cor:contr}
In conditions of Lemma~\ref{lemma:half} we have:
$\ell(g_0)+\ell(g_1) \le \ell(g)+1$.
\end{cor}

The above corollary is not tight and can be improved
in certain cases. The following exercise give bounds
in the other direction, limiting potential extensions
of Corollary~\ref{cor:contr}.

\begin{exc}\label{exc:contr}
In conditions of Lemma~\ref{lemma:half} we have:
$\ell(g) \le 2\ell(g_0)+2\ell(g_1) +50$.
\end{exc}

This result can be used to show that~$\ga_\Gr \suc \exp(\sqrt{n})$.
The proof is more involved that of other exercises; it will not
be used in this paper.

\begin{exc}\label{exc:2^k}
Prove that every element $g \in \Gr$ has order $2^k$, for some
integer~$k$. \ {\rm Hint:} use induction to reduce the problem
to elements $g_0$, $g_1$ (cf. Lemma~\ref{lemma:half}).
\end{exc}

\vskip.8cm

\section{The word problem} \label{sec:word}

The classical \emph{word problem} can be formulated as follows:
given a word $w = s_{i_1} \cdots s_{i_n}$ in generators
$s_{j} \in S$, decide whether
this product is equal to~$\id$ in~$G = \lan S\ran$.
To set up the problem carefully one would have to describe
presentation of the group and allowed operations~\cite{Ha}.
We skip these technicalities in hope that the reader has
an intuitive understanding of the problem.

Now, from the algorithmic point of view the problem is undecidable,
i.e. there is no Turing machine which can resolve it in finite
time for every group.  On the other hand, for certain groups the
problem can be solved very efficiently, in time polynomial in the
length~$n$ of the product.
For example, in the \emph{free group}
$F_k = \lan x_1^{\pm 1},\ldots, x_k^{\pm 1}\ran$
the problem can be solved in linear time: take a product~$w$
and repeatedly cancel every occurrence of~$x_i x_i^{-1}$
and~$x_i^{-1}x_i$, $1 \le i \le k$;
the product~$w=\id$ if and only if the resulting
word is empty.  Since every letter is cancelled at most once
and new letters are not created, the algorithm takes~$O(n)$
cancellations.

\begin{exc} \label{exc:free-linear}
Note that the above algorithm as defined one might need
to search for the next cancellation, increasing the
complexity of the algorithm to as much as $O(n^2)$.
Modify the algorithm by using a single stack to prove
that word problem in~$F_k$ can in fact be solved in
linear time.
\end{exc}

The class of groups where the word problem can be solved in
linear number of cancellations is called \emph{word hyperbolic}.
This class has a simple description and many group theoretic
applications~\cite{Gr}.
The following result shows that word problem can be resolved in~$\Gr$ in
nearly linear time\footnote{In computer science literature \emph{nearly
linear time} usually stands for $O(n \log^kn)$, for some fixed~$k$.}.

\begin{thm} \label{thm:word}
The word problem in~$\Gr$ can be solved in $O(n \log n)$ time.
\end{thm}

\begin{proof} Consider the following algorithm.  First, cancel
products of $b,c,d$ to write the word as
$w = (a) \any a \any \cdots \any a \any (a)$. If the number~$\pi(w)$
of $a$'s is odd, then the product~$w \ne_\Gr \id$.  If the~$\pi(w)$ is
even, use the rewriting rules (proof of Lemma~\ref{lemma:half})
to obtain words~$w_0=\Phi_0(w)$ and~$w_1=\Phi_1(w)$ (which may
no longer by reducible).
Recall that the product $w=_\Gr\id$ if and only if~$w_0, w_1 =_\Gr \id$.
Now repeat the procedure for the words~$w_0,w_1$ to obtain
words~$w_{00},w_{01},w_{10},w_{11}$, etc.
Check that~$w=_\Gr\id$ if and only if eventually
all the obtained words are trivial.

Observe that the length of each word~$w_i$ is at most~$(n+1)/2$.
Iterating this bound, we conclude that the number of
`rounds' in the algorithm of constructing smaller and smaller words
is~$O(\log n)$.  Therefore, each letter is replaced at most~$O(\log n)$
times and thus the algorithm finishes in~$O(n \log n)$ time.
\end{proof}

\begin{Remark} {\rm For every reduced decomposition as above one can
construct a binary tree of nontrivial words $w_{i_1i_2\dots i_r}$.
The distribution of \emph{height} and \emph{shape} (profile) of
these trees is closely connected to the growth function~$\ga_\Gr$.
Exploring this connection is of great interest,
but lies outside the scope of this paper.}
\end{Remark}

\vskip.8cm

\section{Subexponential growth of~$\Gr$} \label{sec:upper}

In this section we prove the second half of Theorem~\ref{thm:main} by
establishing the upper bound on the growth function~$\ga$ of
group~$\Gr$ with generators~$S = \{a,b,c,d\}$.  The proof relies
on the technical Cancellation Lemma which will be stated here and
proved in the next section.

\smallskip

Let $\Hr_3:=\St_\Gr(3)$ be the stabilizer of vertices on the third
level, and recall that the index $[\Gr:\Hr_3] \le 2^7 = 128$
(Exercise~\ref{exc:stab-ind}). There is a natural embedding
$$\psi_3 \, : \, \Hr_3 \, \longrightarrow \, \Gr_{\0\0\0} \times \Gr_{\0\0\1} \times
\ldots \times \Gr_{\1\1\1}$$
(see Section~\ref{sec:inf}).   By self-similarity, each of the
eight groups in the product is isomorphic:
$\Gr_{\mathbf{ijk}} \simeq \Gr$, where~$\mathbf{i},\mathbf{j},\mathbf{k}
\in \{\0,\1\}$. These isomorphisms are obtained by restrictions of
natural maps: $\io_v^{-1}: \Aut(\T_v) \to \Aut(\T)$, where~$v\in \T$.
Now combine~$\psi_3$ with the map
$(\io_{\0\0\0}^{-1},\io_{\0\0\1}^{-1},\ldots,\io_{\1\1\1}^{-1})$ we
obtain a group homomorphism $\chi: \Hr_3 \to \Gr^8$ written as
$\chi(h) = (g_{000},g_{001},\ldots,g_{111})$, where $h \in \Hr_3$
and $g_{ijk} \in \Gr$.

It follows easily from Corollary~\ref{cor:contr}  that
$\ell(g_{000}) + \ell(g_{001}) + \ldots + \ell(g_{111}) \le \ell(h) + 7$.
The following result is an improvement over this bound:

\begin{lemma}[Cancellation Lemma] \label{lemma:can}
Let~$h \in \Hr_3$.  In the notation above we have:
$$\ell(g_{000}) + \ell(g_{001}) + \ldots + \ell(g_{111}) \, \le
\, \frac56 \, \ell(h) + 8.$$
\end{lemma}

\smallskip

We postpone the proof of Cancellation Lemma till next section.
Now we are ready to finish the proof of the Main Theorem.

\begin{prop} Group~$\Gr$ has subexponential growth.  Moreover,
$\ga_\Gr(n) \pre \exp(n^\nu)$ for some~$\nu < 1$.
\end{prop}

\begin{proof}
All elements $g\in \Gr$
can be written as~$g = u \cdot h$ where $h \in \Hr_3$ and $u$~is
a coset representative of~$\Gr/\Hr_3$. Since $[\Gr:\Hr_3]\le 128$,
there are at most~128 such elements~$u$. Note that we can choose
elements~$u$ which have length at most~127
in~$S = \{a,b,c,d\}$,  since all prefixes of a reduced
decomposition can be made to lie in distinct cosets.  The
decomposition~$h = u^{-1}g$ then gives~$\ell(h) \le \ell(g) + 127$.

Now write $g  =  u \ts h = u \ts g_{000}\ts g_{001}\cdots g_{111}$.
The Cancellation Lemma gives:
$$\sum_{ijk} \ell(g_{ijk}) \, \le \, \, \frac56 \, \ell(h) + 8
\, \le  \, \frac56 \, \bigl(\ell(g)+127\bigr) + 8 \, < \,
\frac56 \, \ell(g)+114.$$
Putting all this together we conclude:
$$\ga(n) \, \le \, 128  \sum_{(n_1,\ldots,n_8)} \, \ga(n_1) \cdots \ga(n_8),
$$
where the summation is over all integer $8$-tuples with
$n_1+\ldots+n_8 \le \frac56 n + 114$.
Set~$m=n+137$ so that $\frac56 \ts n + 114 < \frac56 \ts m$.
Now note that
$$\ga(m) \, = \, \ga(n+137) \, \le \, \ga(n) \cdot |S|^{137}
\, = \,4^{137} \, \ga(n).
$$
Therefore, we have:
$$
\ga(m) \, \le \, 4^{137} \, \ga(n) \, \le
\,  4^{137} \cdot 128 \cdot \ga^{\star \,8}\left(\frac56\, n + 114\right)
\, \le \, 2^{281} \, \ga^{\star \,8}\left(\frac56\,m\right).
$$
From here and the Upper Bound (Lemma~\ref{lemma:ub}) we obtain the result.
\end{proof}

\medskip
Recall that subexponential growth of~$\Gr$ is shown in
Corollary~\ref{cor:lower}.
This completes the proof of Theorem~\ref{thm:main}. \ $\sq$

\vskip.8cm

\section{Proof of the Cancellation Lemma}
\label{sec:can-lemma}

Fix a reduced decomposition~$(a) \any a \any a \cdots \any (a)$
of~$h\in \Hr_3$, and denote this decomposition by~$w$.  Apply rewriting
rules $\Phi_0$ and $\Phi_1$ to~$w$ obtain words~$w_0$ and~$w_1$.
At this moment remove all identities~$\id$.  Then apply these rules
again to obtain~$w_{00},w_{01},w_{10}$ and~$w_{11}$, and remove the
identities~$\id$.  Finally, repeat this once again to obtain
words~$w_{000},w_{001},\ldots,w_{111}$.
Following the proof of Theorem~\ref{thm:word}, all these words
give decompositions of elements $g_0,g_1$, then~$g_{00},\ldots,g_{11}$,
and~$g_{ijk} \in \Gr_{\textbf{ijk}}$, respectively.
Note that these decompositions are not necessarily reduced,
so for the record:
$$(\maltese) \ \ \
\ell(g_i) \le |w_i|, \ \ \ell(g_{ij}) \le |w_{ij}|, \ \
\ell(g_{ijk}) \le |w_{ijk}|, \ \ \ \text{for all} \ \ i,j,k\in \{0,1\},
$$
where $|u|$ denotes the length of the word~$u$. Also,
by Corollary~\ref{cor:contr} we have:
$$(\diamondsuit) \quad
\aligned
\ell(g_0)+\ell(g_1) & \le \ell(h)+1, \\
\ell(g_{00}) + \ldots + \ell(g_{11}) & \le \ell(g_0)+\ell(g_1)+2, \\
\ell(g_{000}) + \ell(g_{001}) + \ldots + \ell(g_{111}) &
\le\ell(g_{00}) + \ldots + \ell(g_{11})+4.
\endaligned \quad
$$
To simplify the notation, consider the following
 concatenations of these words:
$$w'= w_0\cdot w_1,  \ \ w'' = w_{00}\cdots w_{11}, \ \ \,
\text{and} \ \
w'''= w_{000}\cdot w_{001} \cdots w_{111}.
$$
By construction of the rewriting rules, since the only possible
cancellation happens when~$d \to \id$ we have:
$|w^\pr|  \le |w| +1 -|w|_d$, where~$|w|_d$ is the number
of letters~$d$ in~$w$. Indeed, simply note that each letter~$d$
in~$w$ is cancelled by either~$\Phi_0$ or~$\Phi_1$.
Unfortunately we cannot iterate this
inequality as the words~$w_i$ are not reduced. Note on
the other hand, that each letter~$c$ in~$w$ produces one
letter~$d$ in~$w'$ and each of the latter is cancelled again
by either~$\Phi_0$ or~$\Phi_1$.  Finally, each letter~$b$ in~$w$
produces one letter~$c$ in~$w'$, which in turn produces letter~$d$
in~$w''$, and each of the latter is cancelled again
by either~$\Phi_0$ or~$\Phi_1$.  Taking into account the
types of decompositions we obtain:
$$(\heartsuit) \ \qquad
\aligned
|w'| \, & \le \, |w| +1 -|w|_d\,, \\
|w''| \, & \le \, |w| +3 -|w|_c\,,\\
|w^{\pr\pr\pr}| \, & \le \, |w| +7 -|w|_b\,.
\endaligned \quad
$$
Since $|w|_b+|w|_c+|w|_d \ge (|w|-1)/2$, at least one of the
numbers $|w|_\any > |w|/6-1$. Combining this
with~$(\heartsuit)$,~$(\diamondsuit)$, and~(\maltese) we conclude:
$$\aligned
&\ell(g_{000}) + \ell(g_{001}) + \ldots + \ell(g_{111}) \,
\le \, \max\bigl\{|w'|+2+4, |w''|+4,|w'''|\bigr\}  \, \\
& \ \ \  \le \,
|w| + 7 - \max_{\any \in \{b,c,d\}} |w|_\any \, \le \,
|w| + 7 - (|w|/6-1) \, = \, \frac56 \, \ell(h) + 8,
\endaligned
$$
as desired. \ $\sq$


\vskip.8cm

\section{Further developments, conjectures and open problems}
\label{sec:open}

Much about groups of intermediate growth remains open.
Below we include only the most interesting results and
conjectures which are closely connected to the material
presented in this paper. Everywhere below we refer to
surveys~\cite{BGN,BGS,GNS} and books~\cite{Ha,Ne} for
details and further references.

\smallskip

Let us start by saying that the Upper Bound and Lower Bound lemmas
can be used to obtain effective bounds on the growth function
of~$\Gr$. Although considerably sharper bounds are known, the
exact asymptotic behavior of~$\ga_\Gr$ remain an open problem.
Unfortunately, we do not even know whether it makes sense to say
that~$\ga_\Gr$ has growth~$\exp(n^\al)$ for some fixed~$\al>0$:

\begin{conj} \label{conj:irate}
Let $\ga=\ga_\Gr$ be the growth function of group~$\Gr$.
Prove that there exists a limit~$\,
\al = \lim_{n\to \infty} \log_n \log \ts \ga(n)$.
\end{conj}

In fact, the limit as in the conjecture is not known to exist
and satisfy $0 <\al < 1$ for \emph{any} finitely generated group.
Alternatively, the extend to which results for~$\Gr$ generalize
to other groups of intermediate growth remains unclear as well.
Although there are now constructions of groups with
subexponential growth function~$\ga(n) \sim e^{n \ts (1-o(1))}$,
there are no known examples of groups with
superpolynomial growth function~$\ga(n) \sim \exp(n^{o(1)})$.
The following result has been established for a large
classes of groups, but not in general:


\begin{conj} \label{conj:lower}
Let~$G$ be a group of intermediate growth, and
let~$\ga_G(n)$ be its growth function.  Then
$\ga_S(n) \suc \exp(n^\al)$ for some~$\al>0$.
\end{conj}

Moving away from the bounds on the growth, let us mention
that group~$\Gr$ is not finitely presented.  Existence of
finitely presented groups of intermediate growth is a major
open problem in the field, and the answer is believed to be
negative:

\begin{conj} \label{conj:fpres}
There are no finitely presented groups of intermediate growth.
\end{conj}

Our final conjecture may seem technical and unmotivated as stated.
If true it resolves positively the ``$p_c<1$'' conjecture of
Benjamini and Schramm on percolation on Cayley graphs.

\begin{conj} \label{conj:perc}
Every group~$G$ of intermediate growth contains two
infinite subgroups~$H_1$ and~$H_2$ which commute with
each other:~$[h_1,h_2]=\id$ for all~$h_1\in H_1$
and~$h_2\in H_2$.
\end{conj}

We refer to~\cite{MP1} for an overview of this conjecture and
its relation to groups of intermediate growth.

\vskip .7cm

{\bf Acknowledgements.} \
We would like to thank Tatiana Nagnibeda and Roman Muchnik
the interest in the subject and engaging discussions.
Both authors were partially supported by the NSF.




\newpage

\section{Appendix} \label{sec:app}

\medskip

{\bf Proof of the Lower Bound Lemma.} \
To simplify the notation, let us extend definition of~$f$ to
the whole line~$f: \rr_+ \to \rr_+$
by setting $f(x):=f(\lfloor x\rfloor)$.
Let $\pi(n) = \ln f(n)$. Clearly, $\pi(n)$ is
monotone increasing, and $\pi(n) \to \infty$ as $n \to \infty$.
By definition, condition $f \suc f^m$ gives
$f(n) \ge C\,f^m(\al\ts n)$ for some $C,\al>0$. Write this
as
$$(\dv) \ \ \ \ \ \
\pi(n) \, \ge \,  c \, + \, m \, \pi(\al \ts n)\,, \ \ \
$$
where $c = \log C$.  Let us first show that $\al <1$. Indeed,
if~$\al \ge 1$, we have:
$$(\dv\dv) \ \ \
m \, \pi(\al\ts n) \, - \, \pi(n) \, > \,
m \, \pi(n) \, - \, \pi(n) \, = \, (m-1) \, \pi(n) \, \to \, \infty \ \
\text{as} \ \ n \to \infty.
$$
On the other hand, $(\dv)$ implies that the l.h.s. of $(\dv\dv)$
is $\le -c$, a contradiction.

Applying $(\dv)$ repeatedly to itself gives us:
$$\aligned
\pi(n) \, & \ge \,  c  +  m \ts \pi(\al\ts n)
\, \ge \, c  +  m \ts (c+ m \ts \pi(\al\ts n)) \,  \ge \, \ldots \\
& \ge \, m^k\pi(\al^kn) + c (1+m+\ldots+m^{k-1}).
\endaligned
$$

Suppose~$c\ge 0$.  Take $k=\lfloor (\log n - 1)/\log1/\al\rfloor$.
Then $\pi(\al^kn) \ge \log(\al^kn) \ge 1$ and
$\pi(n) \, \ge \, m^k \, \ge A \ts n^\nu$, where
$m^{1/\log\al} \ge A \ge m^{(1/\log\al)-1}$ and
$\nu = (\log m)/(\log 1/\al)>0$.

Suppose now~$c< 0$ and recall that~$m\ge 2$. Then
$\pi(n) \ge m^k\bigl(\pi(\al^kn) + c\bigr)$.
Take $k=\lfloor (\log n +c- 1)/\log1/\al\rfloor$.
Then $\pi(\al^kn) +c \ge \log(\al^kn) +c \ge 1$, and
$\pi(n) \, \ge \, m^k \, \ge A \ts n^\nu$, where
$m^{-(c-1)/\log\al} \ge A \ge m^{-1-(c-1)/\log\al}$
and~$\nu$ as above.

Therefore, in both cases we have
$f(n) = \exp\ts\pi(n) \ge \exp(A\ts n^\nu)$
for some~$A,\nu >0$, as desired.
\ $\sq$

\vskip.6cm

{\bf Proof of the Upper Bound Lemma.} \
We prove the result by induction on~$n$. Suppose
$\pi(n) := \log f(m) \le A\ts n^\nu$.  We have:
$$
f(n) \, \le \, C\ts f^{\star \, k}(\al \ts n) \, = \, C
\sum_{(n_1,\ldots,n_k)} \, f(n_1)\cdots f(n_k),
$$
where the summation is over all~$n_1+\ldots+n_k \le \al \ts\ts n$.
Clearly, the number of terms of the summation is at most~$(\al \ts n)^k$.
Also, for each product in the summation we have by induction:
$$\aligned
& \log\bigl(f(n_1)\cdots f(n_k)\bigr) \,  \le \,
\pi(n_1) + \ldots + \pi(n_k) \, \le
\, A \ts (n_1^\nu + \ldots + n_k^\nu)\\
 & \qquad \quad \le \, A \ts k\ts (\al \ts n/k)^{\nu} \, \le \,
A \ts n^\nu \cdot \left[k\,\left(\frac{\al}{k}\right)^\nu\right] \, = \,
A \ts n^\nu \cdot (1 - \ve),
\endaligned
$$
where~$\ve = \ve(k,\al)>0$, and where~$\nu<1$ is large enough.
Therefore,
$$(\diamond) \ \ \ \
\aligned
\pi(n) \, & = \, \log f(n) \, \le \, \log C + \log (\al \ts n)^k +
A \ts n^\nu \cdot (1 - \ve) \\
& \le \,(\log C + k\log \al + k \log n) +
A \ts n^\nu \cdot (1 - \ve) \, \le \, A\ts n^\nu
\endaligned
$$
for~$A$ large enough.  Setting~$A$ large enough to satisfy~$(\diamond)$
and the base of induction, we obtain the result.  \ $\sq$

\end{document}